\documentclass[11pt,a4paper]{article}

\usepackage{amsmath,amssymb,latexsym}
\usepackage{graphicx}
\usepackage{epsfig}
\usepackage{longtable}
\usepackage{amsthm}
\usepackage{xcolor}

\newtheorem{theorem}{Theorem}
\newtheorem{lemma}[theorem]{Lemma}

\newtheorem{definition}{Definition}

\newtheorem{remark}{Remark}
\newtheorem{problem}{Open Problem}

\begin{document}

\begin{center}
{\Large\bf Supports of quasi-copulas}
\end{center}

\begin{center}
{\large\bf Juan Fern\'andez-S\'anchez$^{\rm a,b}$, Jos\'e Juan Quesada-Molina$^{\rm c}$ and Manuel \'Ubeda-Flores$^{\rm d,}$}\footnote{Corresponding author}
\end{center}

\begin{center}
{\it $^{\rm a}$Research Group of Theory of Copulas and Applications, University of Almer\'{\i}a,
04120 Almer\'{\i}a, Spain.\\
$^{\rm b}$Interdisciplinary Mathematics Institute, Complutense University of Madrid, 28040 Madrid, Spain.\\
 {\rm juanfernandez@ual.es}\\
$^{\rm c}$Department of Applied Mathematics, University of Granada, 18071 Granada, Spain.\\{\rm jquesada@ugr.es}\\
$^{\rm d}$Department of Mathematics, University of Almer\'{\i}a, 04120 Almer\'{\i}a, Spain.\\
{\rm mubeda@ual.es}}
\end{center}
\smallskip

\begin{center}{\bf Abstract}\end{center}

\noindent
It is known that for every $s\in]1,2[$ there is a copula whose support is a self-similar fractal set with Hausdorff ---and box-counting--- dimension $s$. In this paper we provide similar results for (proper) quasi-copulas, in both the bivariate and multivariate cases.
\bigskip

\medskip


\section{Introduction}

$\phantom{999}$Aggregation of pieces of information coming from different sources is an important task in expert and decision support systems, multi-criteria decision making, and group decision making. Aggregation operators are precisely the mathematical objects that allow this type of information fusion; and well-known operations in logic, probability theory, and statistics fit into this concept (for an overview, see \cite{Beliakov2007,Grabish2009}). Conjunctive aggregation operators, i.e., those aggregation operators that are bounded by the minimum, constitute an
important class of operators that includes copulas and quasi-copulas.

(Bivariate) quasi-copulas were introduced in the field of probability (see \cite{Alsina93}) and were characterized in \cite{Genest99}. They are also used in aggregation processes because they ensure that the aggregation is stable, in the sense that small error inputs correspond to small error outputs. In the last few years these functions have attracted an increasing interest by researchers in some topics of fuzzy sets theory, such as preference modeling, similarities and fuzzy logics. For a complete overview of quasi-copulas, we refer to \cite{Arias2020,Sempi2017}.

Copulas, (bivariate) probability distribution functions with uniform margins on $[0,1]$ restricted to the unit square, are a subclass of quasi-copulas. The importance of copulas in probability and statistics comes from \textit{Sklar's theorem} \cite{Sklar59}, which states that the joint distribution $H$ of a pair of random variables $(U,V)$ and the corresponding (univariate)
marginal distributions $F$ and $G$ are linked by a copula $C$ in the following manner:
\begin{equation*}
H(x,y)=C\left(F(x),G(y)\right)\,\, \mathrm{for}\,\,\mathrm{all}\,\,\left(x,y\right)\in[-\infty,\infty]^2.
\end{equation*}
If $F$ and $G$ are continuous, then the copula is unique; otherwise, the copula is uniquely determined on Range\,$F\times$Range\,$G$ (\cite{deAmo2012}). For a comprehensive review on copulas, we refer to the monographs \cite{Durante2016book,Nelsen2006}.

The construction of (bivariate) proper quasi-copulas, i.e., quasi-copulas that are not copulas, from the so-called proper quasi-transformation square matrices shows that these functions do not induce a stochastic signed measure on $[0,1]^2$ (see \cite{FerRodUb2011}, and the case of higher dimensions in \cite{FerUb2014,Nelsen2010}).

A \textit{fractal} is a set whose topological dimension is less than its Hausdorff dimension. In \cite{Fredricks2005}, copulas whose supports are fractals were constructed in such a way that for every $s\in ]1,2[$ it is possible to find a copula whose support is a self-similar fractal set with Hausdorff ---and box-counting--- dimension $s$ (see \cite{Trutschnig2012} for a generalization to the multivariate case and \cite{deAmo2012b} for other examples).

Our purpose in this paper is to generate quasi-copulas whose support ---an extension of the concept of support of a copula--- is a self-similar fractal set and has Hausdorff dimension $s$.

The paper is organized as follows. After some preliminaries concerning copulas, quasi-copulas and signed measures (Section \ref{prelim}), in Section \ref{qtransform} we review some concepts related to quasi-transformation matrices. In Section \ref{fractals} we introduce the concept of support of a quasi-copula in order to study quasi-copulas whose support is self-similar and has Hausdorff dimension in $]1,2]$. In Section \ref{multi} we tackle the multidimensional case. Finally, Section \ref{conclus} is devoted to conclusions.

\section{Preliminaries}\label{prelim}

$\phantom{999}$A (bivariate) \textit{copula} is a function $C\colon [0,1]^{2}\longrightarrow [0,1]$ which satisfies:

\begin{itemize}
\item[(C1)] the boundary conditions $C(t,0)=C(0,t)=0$ and $C(t,1)=C(1,t)=t$ for all $t$ in $[0,1]$, and

\item[(C2)] the \textit{2-increasing property}, i.e., $V_{C}([u_{1},u_{2}]\times[v_{1},v_{2}])=C(u_{2},v_{2})-C(u_{2},v_{1})-C(u_{1},v_{2})+C(u_{1},v_{1})\ge 0$ for all $u_{1},u_{2},v_{1},v_{2}$ in $[0,1]$ such that $u_{1}\le u_{2}$ and $v_{1}\le v_{2}$.
\end{itemize}

$V_C(R)$ is usually called as the $C$-\textit{volume} of the rectangle $R$; and in the sequel we also consider the $C$-volume of a rectangle for real-valued functions on $[0,1]^2$ which may not be copulas. The set of copulas is denoted by $\mathcal{C}$.

A \textit{quasi-copula} is a function $Q\colon[0,1]^{2}\longrightarrow [0,1]$ which satisfies condition (C1) of copulas, but in place of (C2), the following weaker conditions:

\begin{itemize}
\item[(Q2)] $Q$ is non-decreasing in each variable, and

\item[(Q3)] the \textit{Lipschitz} condition $|Q(u_{1},v_{1})-Q(u_{2},v_{2})|\le |u_{1}-u_{2}|+|v_{1}-v_{2}|$ for all $u_{1},v_{1},u_{2},v_{2}$ in $[0,1]^{2}$
\end{itemize}
(see \cite{Genest99}).

Let ${\mathcal{Q}}$ denote the set of quasi-copulas. While every copula is a quasi-copula, there exist \textit{proper} quasi-copulas, i.e., quasi-copulas which are not copulas.

For each $n=1,2$, let $\mathfrak{B}^n$ denote the Borel $\sigma$-algebra for $[0,1]^n$, and let ${\mathcal{S}}^n$ denote the measurable space $\left([0,1]^n,\mathfrak{B}^n\right)$. It is known that every bivariate copula $C$ induces a probability measure $\mu_C$ on ${\mathcal{S}}^2$ such that $\mu_C\left([0,1]\times A\right)=\mu_C\left(A\times [0,1]\right)=\lambda(A)$ for every set $A$ in $\mathfrak{B}^1$, where $\lambda$ denote the Lebesgue measure in $\mathbb{R}$; i.e., $\mu_C$ is a \emph{doubly stochastic measure} (see, e.g., \cite{Nelsen2006}). This measure is characterized by the fact that $\mu_C(R)=V_C(R)$ for every rectangle $R\subseteq[0,1]^2$. It also satisfies $0\le \mu_C(D)\le 1$ for every set $D$ in $\mathfrak{B}^2$. The \textit{support} of a copula $C$
---denoted by supp$(C)$ or supp$\left(\mu_C\right)$--- is the complement of the union of all open subsets of $[0,1]^2$ with $\mu_C$-measure zero. When we refer to ``mass'' on a set, we mean the value of $\mu_C$ for that set.

A \emph{signed measure} $\mu$ on ${\mathcal{S}}^2$ is an extended real valued, countably additive set function on $\mathfrak{B}^2$, such that $\mu(\emptyset)=0$ and $\mu$ assumes at most one of the values $\infty$ and $-\infty$. Equivalently, $\mu$ is the difference between two (positive) measures $\mu_1$ and $\mu_2$ on ${\mathcal{S}}^2$, such that at least one of them is finite. \textit{Hahn's decomposition theorem} guarantees that we can take $\mu_1\bot \mu_2$, i.e., $\mu_1$ and $\mu_2$ are \textit{mutually singular}; in other words, there exists a measurable set $A$ such that $\mu_1(A)=\mu_2(A^c)=0$. For signed measures, we have $\mathrm{supp}(\mu)=\mathrm{supp}(\mu_{1})\cup \mathrm{supp}(\mu_2)$, with the condition that $\mu_{1}\bot \mu_2$. Just like positive measures, a signed measure $\mu $ on ${\mathcal{S}}^{2}$ is said to be \emph{doubly stochastic} if $\mu([0,1]\times A)=\mu (A\times \lbrack 0,1])=\lambda (A)$ for every set $A$ in
$\mathfrak{B}^{1}$. For more details, see \cite{Halmos1974,Nelsen2010}.

A \emph{hyperbolic IFS} (\textit{iterated function system}) consists of a complete metric space $(X,d)$ together with a finite set of contraction mappings $(w_{n}:X\longrightarrow X)$. The hyperbolic IFS $\{X;(w_{n})\}$ satisfies \emph{Moran's open set condition} if there exists a nonempty open subset $U$ of $X$ for which $w_{n}(U)\cap w_{m}(U)=\emptyset $ whenever $n\neq m$ and $w_{n}(U)\subseteq U$ for all $n$. Given an IFS$\ \{X;(w_{n})\}$ and a set $A\subseteq X$, we define $W(A)=\cup_{n}w_{n}(A)$.

Let $(X, d)$ be a metric space and $d(x, y)$ the distance between $x$ and $y$ in $X$. A map $S \colon X\longrightarrow X$ is called a \textit{similarity with ratio c} if there exists a number $c \in]0,1[$ such that $d(S(x), S(y)) = c\cdot d(x, y)$ for all $x, y\in X$. A set $E$ is said to be \textit{self-similar} if it can be expressed as $E=\cup_{k=1}^{m}S_k(E)$, for some
similarities $S_1,S_2,\ldots,S_m$. See \cite{Hutchinson1981,Mandelbrot1977} for more details.

There are many ways to define fractal dimension and not all definitions are equivalent. Two of the most well-known definitions of fractal dimension are the Hausdorff dimension and the box-counting dimension. We recall both concepts (see \cite{Falconer2014} for further details).

Let $n\ge 2$ be a natural number, and let $\left(\mathbb{R}^{n},d\right)$ be a metric space. For any bounded subset $E\subset \mathbb{R}^{n}$, let $\mathrm{diam}(E)$ denote the diameter of $E$, i.e., $\mathrm{diam}(E)=\sup\{d(x,y) : x,y\in E\}$. For any $\delta\in ]0,\infty]$, let $\left\{E_i\right\}_{i\ge 1}$ be a $\delta$-cover of $E$, i.e., $E\subseteq \cup_{i\ge 1}E_i$, with $\mathrm{diam}\left(E_i\right)\le\delta$ for all $i$. For any $\alpha\in [0,\infty[$ consider the outer measure
\begin{equation*}
{\mathcal{H}}_{\delta}^{\alpha}(E):=\inf\left\{\sum_{i=1}^{\infty}\left(\mathrm{diam}\left(E_i\right)\right)^{\alpha}: \left\{E_i\right\}_{i\ge 1}\,\text{is a}\,\, \delta\text{-cover of}\,\,E\right\}.
\end{equation*}
The $\alpha$-dimensional Hausdorff measure of $E\subset \mathbb{R}^{n}$ is defined as
\begin{equation*}
{\mathcal{H}}^{\alpha}(E):=\lim_{\delta\to 0}{\mathcal{H}}_{\delta}^{\alpha}(E),
\end{equation*}
and the \textit{Hausdorff dimension} of $E$ is the number defined by
\begin{equation*}
\dim_{\mathcal{H}}(E):=\inf\left\{\alpha\ge 0:{\mathcal{H}}^{\alpha}(E)=0\right\}=\sup\left\{\alpha\ge 0:{\mathcal{H}}^{\alpha}(E)=\infty\right\}.
\end{equation*}
Let $N_{\delta}(E)$ be the smallest number of sets of diameter at most $\delta$ that cover $E$. The \textit{box}-\textit{counting dimension} of $E$ is defined by
\begin{equation*}
\dim_{\mathcal{B}}(E):=-\lim_{\delta\to 0}\frac{\ln N_\delta(E)}{\ln \delta}.
\end{equation*}
In general, we have $\dim_{\mathcal{H}}(E)\le \dim_{\mathcal{B}}(E)$. The Hausdorff dimension and the box-counting dimension are sometimes equal ---for example, when we have a self-similar set---; but, for instance, for the set $A=\mathbb{Q} \cap [0, 1]$ we have $\dim_{\mathcal{H}}(A)=0$ and $\dim_{\mathcal{B}}(A)=1$.

The following result, which can be found in \cite{Barnsley1989,Edgard2008}, will be useful for our purposes.

\begin{theorem}
\label{IFS} Let $(X,d)$ be a complete metric space. If $\{X;(w_{n})\}$ is a hyperbolic IFS, then there exists a unique nonempty compact subset $K$ of $X$ for which $K=W(K)$. Moreover, if each $w_n$ is a similarity and Moran's open set condition is satisfied, then $\dim_{\mathcal{H}}(K)=\dim_{\mathcal{B}}(K)=s$ and it is implicitly given by the equation $\sum_{n}c_{n}^{s}=1$,
where, for each $n$, $c_{n}$ is the similarity ratio of $w_{n}$.
\end{theorem}

The set $K$ in Theorem \ref{IFS} is called the \textit{invariant set} of $\{X;(w_{n})\}$ and we denote it by $K_{W}$. Moreover, if $X$ is compact, then it is satisfied $K_W=\cap_{k=0}^{\infty}W^{k}(X)$.

\section{Quasi-transformation matrices}

\label{qtransform}

$\phantom{999}$We recall the following definition in order to represent quasi-copulas in terms of matrices.

\begin{definition}[\protect\cite{FerRodUb2011}]\label{quasi-transform}
A \emph{quasi-transformation matrix} is a square matrix $T=(t_{ij})$, $i,j=1,2,\ldots,m$, with the column index first and the rows ordered from bottom to top, satisfying the following conditions: (a) the sum of all entries in $T$ is $1$; (b) the sum of the entries in any row or column of $T$ is positive; and (c) the sum of the entries in any submatrix $T$ that contains at least one entry from the first or last row or column of $T$, is nonnegative.
\end{definition}

\begin{remark}
We want to stress that the original definition in \cite{FerRodUb2011} includes the extra condition ``every entry $t_{ij}$ is in $[-1/3,1]$,'' but this condition is a consequence of the other three conditions, so we omit it.
\end{remark}

A quasi-transformation matrix $T$ is \emph{proper} if some entry in $T$ is negative, for instance, the matrix $T_0$ given by
\begin{equation}
T_0=
\begin{pmatrix}
0 & 1/3 & 0 \\
1/3 & -1/3 & 1/3 \\
0 & 1/3 & 0%
\end{pmatrix}
\!;  \label{T1/3}
\end{equation}
otherwise, if all entries in a quasi-transformation matrix $T$ are nonnegative, then $T$ is a \emph{transformation matrix} (see \cite{Fredricks2005}).

Let $p_{i}$ (respectively, $q_{j}$) denote the sum of the entries in the first $i$ columns (respectively, first $j$ rows) of $T$, for every $i$ (respectively, $j$) in $\{0,1,\ldots m\}$. So $p_{0}=0$ and $p_{m}=1$ (respectively, $q_{0}=0$ and $q_{m}=1$). For every $i,j=1,2,\ldots ,m$, let $R_{ij}$ be the box in $[0,1]^{2}$ given by $R_{ij}=[p_{i-1},p_{i}]\times\lbrack q_{j-1},q_{j}]$. In what follows, we refer to these boxes as the \emph{rectangles associated} with the quasi-transformation matrix $T$. From the definition of quasi-transformation matrix, all these rectangles are not degenerate (they have positive area).

\begin{figure}[tbh]
\begin{center}
\epsfig{file=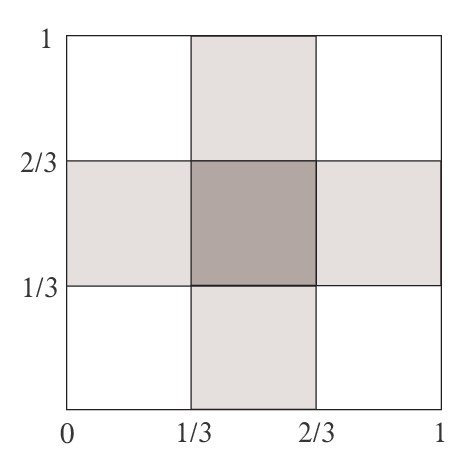,width=5cm}
\end{center}
\caption{The support of $T_{0}(\Pi )$.}
\label{fig:T0}
\end{figure}

For every quasi-transformation matrix $T=(t_{ij})$, with associated rectangles $R_{ij}$, $i,j=1,2,\ldots ,m$, and every quasi-copula $Q$, the \emph{T-transformation} of $Q$, denoted by $T(Q)$, is defined as the function on $[0,1]^{2}$ given by
\begin{equation}
T(Q)(u,v)=\sum_{\displaystyle\overset{i^{\prime }<i}{\scriptstyle j^{\prime}<j}}t_{i^{\prime }j^{\prime }}+\frac{u-p_{i-1}}{p_{i}-p_{i-1}}\sum_{j^{\prime }<j}t_{ij^{\prime }}+\frac{v-q_{j-1}}{q_{j}-q_{j-1}}\sum_{i^{\prime }<i}t_{i^{\prime }j}+t_{ij}\cdot Q\biggl(\frac{u-p_{i-1}}{p_{i}-p_{i-1}},\frac{v-q_{j-1}}{q_{j}-q_{j-1}}\biggr)  \label{TQ}
\end{equation}
for every $i,j=1,2,\ldots ,m$ and every $(u,v)\in R_{ij}$ (the empty sums are considered equal to zero); i.e., for every $i,j=1,2,\ldots ,m$, $T(Q)$ spreads mass $t_{ij}$ on $R_{ij}$ in the same (but re-scaled) way in which $Q$ spreads mass on $[0,1]^{2}$. Figure \ref{fig:T0} shows the support of the proper quasi-copula $T_{0}(\Pi )$, where $\Pi $ is the copula of independent
random variables, i.e., $\Pi (u,v)=uv$ for all $(u,v)\in \lbrack 0,1]^{2}$, where the grey (respectively, dark grey) color corresponds to positive (respectively, negative) mass uniformly (and re-scaled) distributed.

The following theorem is a summary of the main results obtained in the study of this type of quasi-copulas in \cite{FerRodUb2011}.

\begin{theorem}
\label{carac} Let $T$\ be a quasi-transformation matrix of order greater than $1$. Then, there is a unique quasi-copula $Q_{T}$\ for which $T(Q_{T})=Q_{T}$. Moreover, $Q_{T}$ is a proper quasi-copula if, and only if, $T$\ is a proper quasi-transformation matrix. For every proper quasi-transformation matrix $T$, the proper quasi-copula $Q_{T}$ does not induce a doubly stochastic signed measure on $[0,1]^{2}$.
\end{theorem}

\begin{remark}
\label{remQ} Note that, as a trivial consequence of Theorem \ref{carac}, we have that $Q_{T}$ is a copula if, and only if, $T$ is a transformation matrix. Observe also that ${\mathcal{Q}}$ can be divided into two subsets, $\mathcal{Q}_{S}$ and $\mathcal{Q}_{N}$, where $\mathcal{Q}_{S}$ denotes the set of the quasi-copulas that induce a signed measure in $[0,1]^{2}$ and $\mathcal{Q}_{N}:=\mathcal{Q}\backslash\mathcal{Q}_{S}$. Given $Q\in\mathcal{Q}_{S}$, if its induced signed measure is $\mu$, the support of $\mu$ will be called the support of $Q$ and denoted by supp$(Q)$.
\end{remark}

Let $T=(t_{ij})$ be a quasi-transformation matrix of order $m\geq 2$, and let $l\in \mathbb{N}$. Let $R_{ij}=[p_{i-1},p_{i}]\times \lbrack q_{j-1},q_{j}]$, $i,j=1,2,\ldots ,m$, be the rectangles associated with $T$. Given any quasi-copula $Q$, the definition of the quasi-copula $T^{l}(Q)$ ---i.e., $T^{0}(Q)=Q$ and $T^{l}(Q)=T(T^{l-1}(Q))$ for every $l\in \mathbb{N}$--- is associated with $m^{2l}$ rectangles; specifically, they are of the form
\begin{equation*}
R_{i_{1},i_{2},\ldots ,i_{l};j_{1},j_{2},\ldots ,j_{l}}=\omega_{i_{1}j_{1}}\circ \omega _{i_{2}j_{2}}\circ \cdots \circ \omega_{i_{l}j_{l}}\left( [0,1]^{2}\right)
\end{equation*}%
in $[0,1]^{2}$ (for $l=1$ we obtain again the rectangles $R_{ij}$), where
\begin{equation}
\omega _{ij}\colon \lbrack 0,1]^{2}\longrightarrow \lbrack 0,1]^{2},\quad\omega _{ij}(u,v)=\left( \alpha _{i}(u),\beta _{j}(v)\right)  \label{wij}
\end{equation}%
with
\begin{eqnarray*}
\alpha _{i}\!\!\!\!\! &\colon &\!\!\!\![0,1]\longrightarrow \lbrack p_{i-1},p_{i}],\quad \alpha _{i}(u)=p_{i-1}+(p_{i}-p_{i-1})u, \\
\beta _{j}\!\!\!\!\! &\colon &\!\!\!\![0,1]\longrightarrow \lbrack q_{j-1},q_{j}],\quad \beta _{j}(v)=q_{j-1}+(q_{j}-q_{j-1})v; \\
&&
\end{eqnarray*}%
and
\begin{equation*}
V_{Q_{T}}\left( R_{i_{1},i_{2},\ldots ,i_{l};j_{1},j_{2},\ldots,j_{l}}\right) =\prod_{k=1}^{l}t_{i_{k}j_{k}},
\end{equation*}%
where
\begin{equation*}
Q_{T}=\lim_{l\rightarrow \infty }T^{l}(Q)
\end{equation*}%
(see \cite{FerRodUb2011} for details). For the sake of simplicity, we will write $R_{\mathbf{i}_{l}\mathbf{j}_{l}}$ instead of $R_{i_{1},i_{2},\ldots,i_{l};j_{1},j_{2},\ldots ,j_{l}}$, where $\mathbf{i}_{l}=\left(i_{1},i_{2},\ldots ,i_{l}\right) $ and $\mathbf{j}_{l}=\left(j_{1},j_{2},\ldots ,j_{l}\right) $. Moreover, taking into account the definition of the functions $\alpha _{i}$ and $\beta _{j}$, it turns out that the lengths of the sides of $R_{\mathbf{i}_{l}\mathbf{j}_{l}}$ are $\prod_{k=1}^{l}\left( p_{i_{k}}-p_{i_{k-1}}\right) $ and $\prod_{k=1}^{l}\left( q_{i_{k}}-q_{i_{k-1}}\right) $ in the corresponding directions; in particular, these dimensions tend to zero exponentially.

Finally, we provide the following result.

\begin{lemma}\label{lema}
Given a quasi-transformation matrix $T=\left( t_{ij}\right) $, if $t_{i_{l}j_{l}}=0$, then the restriction of $Q_{T}$ to $R_{\mathbf{i}_{l}\mathbf{j}_{l}}$ is given by
\begin{equation}
Q_{T}(u,v)=\gamma _{1}+\gamma _{2}u+\gamma _{3}v,  \label{lineal}
\end{equation}
where $\gamma _{1},\gamma _{2},\gamma_{3}$ are real numbers which depend on $\mathbf{i}_{l},\mathbf{j}_{l}$ and $T $.
\end{lemma}

\begin{proof}
If $(u,v)\in R_{ij}$, Equation (\ref{TQ}) for $Q_{T}$ can be written in the form
\begin{equation*}
Q_{T}(u,v)=\sum_{\displaystyle\overset{i^{\prime }<i}{\scriptstyle j^{\prime}<j}}t_{i^{\prime }j^{\prime }}+\alpha _{i}^{-1}(u)\sum_{j^{\prime}<j}t_{ij^{\prime }}+\beta _{j}^{-1}(v)\sum_{i^{\prime }<i}t_{i^{\prime}j}+t_{ij}\cdot Q_{T}\biggl(\alpha _{i}^{-1}(u),\beta _{j}^{-1}(v)\biggr).
\end{equation*}%
Clearly, if $t_{ij}=0$ in $R_{ij}$, the quasi-copula $Q_{T}$ is given by (\ref{lineal}), for appropriate values $\gamma _{1},\gamma _{2},\gamma _{3}$. Assume the result is true for $l-1$, i.e., if we have a rectangle $R_{\mathbf{i}_{l-1}\mathbf{j}_{l-1}}$ with $t_{i_{l-1},j_{l-1}}=0$, then $Q_{T}$ is given by (\ref{lineal}) in $R_{\mathbf{i}_{l-1}\mathbf{j}_{l-1}}$. Since $w_{i_{1}j_{1}}^{-1}\left( R_{\mathbf{i}_{l}\mathbf{j}_{l}}\right)=R_{i_{2},\ldots ,i_{l};j_{2},\ldots ,j_{l}}$ and $\biggl(\alpha_{i}^{-1}(u),\beta _{j}^{-1}(v)\biggr)=w_{ij}^{-1}(u,v)$, if $(u,v)\in R_{\mathbf{i}_{l}\mathbf{j}_{l}}$, then we have
\begin{equation*}
Q_{T}(u,v)=\sum_{\displaystyle\overset{i^{\prime }<i_{1}}{\scriptstyle j^{\prime }<j_{1}}}t_{i^{\prime }j^{\prime }}+\alpha_{i_{1}}^{-1}(u)\sum_{j^{\prime }<j_{1}}t_{i_{1}j^{\prime }}+\beta
_{j_{1}}^{-1}(v)\sum_{i^{\prime }<i_{1}}t_{i^{\prime}j_{1}}+t_{i_{1}j_{1}}\cdot Q_{T}\biggl(\alpha _{i_1}^{-1}(u),\beta_{j_1}^{-1}(v)\biggr).
\end{equation*}%
Finally, the functions $\alpha _{i_{1}}^{-1}(u)$ and $\beta _{j_{1}}^{-1}(v)$ are affine and $Q_{T}$ is given by (\ref{lineal}) in $w_{i_{1}j_{1}}^{-1}\left(R_{\mathbf{i}_{l}\mathbf{j}_{l}}\right) =R_{i_{2},\ldots,i_{l};j_{2},\ldots ,j_{l}}$ since it is a rectangle with $l-1$ indices and $t_{i_{l}j_{l}}=0$, which completes the proof.
\end{proof}

\section{Quasi-copulas and fractals}\label{fractals}

$\phantom{999}$Let $\mathcal{R}$ denote the set of rectangles of the form $I_1\times I_2$, where $I_1$ and $I_2$ are two open intervals in $[0,1]$.

\begin{definition}\label{openrec}
Let $Q$ be a (proper) quasi-copula. A rectangle $R\in \mathcal{R}$ is \textit{open} $Q$-\textit{null} if $V_{Q}\left( R_{1}\right) =0$ for every rectangle $R_{1}=\left[ a_{1},b_{1}\right] \times \left[ c_{1},d_{1}\right] \subset R$. The \textit{support} of $Q$, denoted by supp$(Q)$, is the complement in $[0,1]^{2}$ of the union of all open $Q$-null rectangles.
\end{definition}

\begin{remark}
We note that the open rectangles in Definition \ref{openrec} are referred to the relative topology (or induced topology); thus, for example, the rectangle $[0,a[\times ]b,c[$ is open. Furthermore, for copulas, this new definition coincides with the classical support of the corresponding doubly stochastic measure.
\end{remark}

As we noticed in the Introduction, in \cite{Fredricks2005} it was proved that there is a copula whose support is a self-similar fractal set with Hausdorff ---and box-counting--- dimension $s$ for every $s\in]1,2[$. We wonder whether this result remains true for the case of quasi-copulas. In this section we answer (affirmatively) this question for the set ${\mathcal{Q}}_N$ and provide a result referring to the Hausdorff dimension and the box-counting dimension in the case of ${\mathcal{Q}}_S$ ---recall Remark \ref{remQ}.

\subsection{Self-similarity, Hausdorff and box-counting dimensions in $]1,2[$}

$\phantom{999}$First, we study the set $\mathcal{Q}_{S}$.

\begin{theorem}
\label{th1} For each $s\in ]1,2[$, there exists $Q\in \mathcal{Q}_{S}$ such that $\dim_{\mathcal{H}}(\mathrm{supp}(Q))=\dim_{\mathcal{B}}(\mathrm{supp}(Q))=s$.
\end{theorem}

\begin{proof}
Let $C\in{\mathcal{C}}$ such that $\dim_{\mathcal{H}}(\mathrm{supp}(C))=\dim_{\mathcal{B}}(\mathrm{supp}(C))=s$, and let $T_0$ be the proper quasi-transformation matrix given by \eqref{T1/3}. Then $Q:=T_0(C)$ is a proper quasi-copula which induces a signed measure formed by five replicas of $\mu _{C}$: four positive replicas with a mass of $1/3$, and a negative one with a mass of $-1/3$. Thus, we have $\dim_{\mathcal{H}}\left(\mathrm{supp}\left(Q\right)\right)=\dim_{\mathcal{B}}\left(\mathrm{supp}\left(Q\right)\right)=s$ since $\mathrm{supp}\left(Q\right)$ is the union
of five self-similar sets to supp$(C)$.
\end{proof}

\begin{problem}
It is an open problem to establish whether it is possible to find a quasi-copula with a self-similar support that induces a signed measure.
\end{problem}

Now we show that the proper quasi-copulas $Q_{T}$, where $T$ is a proper quasi-transformation matrix, satisfy an analogous property to that of Theorem \ref{th1}. The next result provides the support of $Q_{T}$.

\begin{theorem}
\label{lsupp}Given a proper quasi-transformation matrix $T$, the support of the quasi-copula $Q_{T}$ is the invariant set of $\left\{ [0,1]^{2};\left(\omega _{ij}\right) \right\} $, where $\omega _{ij}$ are the functions given by \eqref{wij} and such that $t_{ij}\neq 0$.
\end{theorem}

\begin{proof}
First, we show that the points contained in the invariant set $K_{W}$ of $\left\{ [0,1]^{2};\left( \omega_{ij}\right) \right\}$ belong to supp$(Q_{T}) $. Let $Z$ denote the set of pairs $(i,j)$ such that $t_{ij}\neq 0$, and let $W$ be the function defined in the subsets $A\subseteq \lbrack0,1]^{2}$ given by $W(A)=\cup _{(i,j)\in Z}\omega_{ij}(A)$. Since $K_{W}=\cap _{k}W^{k}([0,1]^{2})$, when $(u,v)\in K_{W}$ we have that for every natural number $l$ there exist $\mathbf{i}_{l}$ and $\mathbf{j}_{l}$ such that $(u,v)\in R_{\mathbf{i}_{l}\mathbf{j}_{l}}$. Since the length of the sides of these rectangles tend to $0$, if $R$ is an open rectangle in $[0,1]^{2}$ such that $(u,v)\in R$, we can take a big enough value $l$ and a rectangle $R_{\mathbf{i}_{l}\mathbf{j}_{l}}$ such that $(u,v)\in R_{\mathbf{i}_{l}\mathbf{j}_{l}}\subset R$. Since $V_{Q_{T}}\left( R_{\mathbf{i}_{l}\mathbf{j}_{l}}\right) =\prod_{k=1}^{l}t_{i_{k}j_{k}}\neq 0$ we have that $R$ is not open $Q_{T}$-null, i.e., $(u,v)\in\text{supp}(Q_{T})$, and thus $K_{W}\subseteq \text{supp}(Q_{T})$.

Now we check that supp$(Q_{T})\subseteq K_{W}$. Assume $(u,v)\notin K_{W}$. We have two possibilities:

\begin{enumerate}
\item Suppose there exists a value $l$ and two vectors $\mathbf{i}_{l}$ and $\mathbf{j}_{l}$ such that $(u,v)$ is in the interior of a rectangle $R_{\mathbf{i}_{l}\mathbf{j}_{l}}$ ---i.e., $(u,v)\in \mathrm{int}\left( R_{\mathbf{i}_{l}\mathbf{j}_{l}}\right) $---, with $t_{i_{l}j_{l}}=0$, then the expression of $Q_{T}(u,v)$ ---which is independent of $Q$--- is of the form (%
\ref{lineal}). If $R=[a_{1},b_{1}]\times \lbrack a_{2},b_{2}]\subset \mathrm{int}\left( R_{\mathbf{i}_{l}\mathbf{j}_{l}}\right) $, we easily obtain $V_{Q_{T}}(R)=0$. Therefore, $\mathrm{int}\left( R_{\mathbf{i}_{l}\mathbf{j}_{l}}\right) $ is an open $Q_T$-null rectangle, whence $(u,v)\notin \mathrm{supp}\left(Q_T\right)$, i.e., supp$\left(Q_{T}\right)\subseteq K_{W}$.

\item If for every $l$ and for every pair of vectors $\mathbf{i}_{l}$ and $\mathbf{j}_{l}$ we have $(u,v)\notin \mathrm{int}\left( R_{\mathbf{i}_{l}\mathbf{j}_{l}}\right) $, then there exists $\varepsilon >0$ small enough such that $]u-\varepsilon ,u+\varepsilon \lbrack \times]v-\varepsilon,v+\varepsilon \lbrack \cap K_{W}=\emptyset $ and a finite union of rectangles ---at most, four---, indexed in $P$, of the form $R_{\mathbf{i}^{p}\mathbf{j}^{p}}$ with $p\in P$ such that $t_{i_{l_{p}}j_{lp}}=0$, satisfying that $\left( u,v\right) $ belongs to
the boundary of each of the regions $R_{\mathbf{i}^{p}\mathbf{j}^{p}}$, and with $]u-\varepsilon ,u+\varepsilon \lbrack \times ]v-\varepsilon ,v+\varepsilon \lbrack $ being in the union of those rectangles. Consider the square $R_{1}:=]u-\varepsilon /2,u+\varepsilon /2[\times ]v-\varepsilon/2,v+\varepsilon /2[$. Given $\left[ a_{1},b_{1}\right] \times \left[a_{2},b_{2}\right] \subset R_{1}$, we apply Equation (\ref{lineal}) to each of the rectangles $\left[a_{1},b_{1}\right] \times \left[ a_{2},b_{2}\right]\cap R_{\mathbf{i}^{p}\mathbf{j}^{p}}$, with $p\in P$. Then we have $V_{Q_{T}}\left( \left[a_{1},b_{1}\right] \times \left[ a_{2},b_{2}\right]\cap R_{\mathbf{i}^{p}\mathbf{j}^{p}}\right) =0$, whence $V_{Q_{T}}\left( \left[ a_{1},b_{1}\right] \times \left[ a_{2},b_{2}\right] \right) =0$, i.e., $R_{1}$ is $Q_{T}$-null and $(u,v)\notin \mathrm{supp}(Q_{T})$, and hence supp$(Q_{T})\subseteq K_{W}$.
\end{enumerate}

We conclude that supp$(Q_{T})=K_{W}$, and this completes the proof.
\end{proof}

\begin{remark}
We note that it will frequently happen that $\omega_{ij}$ are not contraction mappings. Since the only $\omega_{ij}$ to take into account are those for which $t_{ij}\neq 0$, the necessary and sufficient condition for supp$(Q_{T})$ to be a self-similar set is that $t_{ij}\sum_{s=1}^{m}t_{is}=t_{ij}\sum_{s=1}^{m}t_{sj}$.
\end{remark}

The following theorem provides proper quasi-transformation matrices which are not associated with doubly stochastic signed measures on $[0,1]^{2}$.

\begin{theorem}\label{frac1}
Given the proper quasi-transformation matrix
\begin{equation*}
T_{r}=%
\begin{pmatrix}
0 & \frac{r}{15} & \frac{r}{5} & \frac{r}{15} \\
0 & \frac{r}{5} & -\frac{r}{15} & \frac{r}{5} \\
0 & \frac{r}{15} & \frac{r}{5} & \frac{r}{15} \\
1-r & 0 & 0 & 0
\end{pmatrix},
\end{equation*}
with $0<r<1$, the proper quasi-copula $Q_{T_{r}}$ is not associated with a doubly stochastic signed measure, its support is a self-similar set and for every $s\in \left] 1,2\right[ $ there exists a unique $r\in \left] 0,1\right[$ such that $\dim _{\mathcal{H}}\left( \mathrm{supp}\left( Q_{T_{r}}\right)\right) =\dim _{\mathcal{B}}\left( \mathrm{supp}\left( Q_{T_{r}}\right)\right) =s$.
\end{theorem}

\begin{proof}
From Theorems \ref{IFS} and \ref{lsupp} we have $\dim _{\mathcal{H}}\left( \mathrm{supp}\left( Q_{T_{r}}\right) \right) =s(r)$, with $s(r)$ being the solution of the equation $g_{s}(r)=1$, where $g_{s}(r)=\left( 1-r\right) ^{s}+9\left(\frac{r}{3}\right) ^{s}$. Then
\begin{equation*}
g_{s}^{\prime }(r)=s\left[ -\left( 1-r\right) ^{s-1}+3\left( \frac{r}{3}\right) ^{s-1}\right] ,
\end{equation*}
so that $g_{s}^{\prime }(r)<0$ for $r\in ]0,r_{s}[$ and $g_{s}^{\prime}(r)>0 $ for $r\in ]r_{s},1[$ (whence there exists a unique value $r_s\in ]0,1[$ such that $g_{s}^{\prime }\left(r_s\right)=0$), where
\begin{equation*}
r_{s}=\frac{3^{-\frac{1}{s-1}+1}}{1+3^{-\frac{1}{s-1}+1}}.
\end{equation*}
Since for every $s\in \left] 1,2\right[ $ we have
\begin{equation*}
\lim_{r\rightarrow 0^{+}}g_{s}(r)=1<3^{2-s}=\lim_{r\rightarrow
1^{-}}g_{s}(r),
\end{equation*}
then there exists $r\in \left] 0,1\right[ $ such that $g_{s}(r)=1$.

Now, consider the function of two variables
\begin{equation*}
f(r,s)=\left( 1-r\right) ^{s}+9\left( \frac{r}{3}\right) ^{s}.
\end{equation*}
Then, for each $(r,s)\in \left] 0,1\right[ \times \left] 1,2\right[ $, we have
\begin{equation*}
\frac{\partial f}{\partial s}(r,s)=(1-r)^{s}\ln (1-r)+9\left( \frac{r}{3}\right) ^{s}\ln \left( \frac{r}{3}\right) <0;
\end{equation*}
thus, if we implicitly derive $s(r)$ we obtain
\begin{equation*}
s^{\prime }(r)=-\frac{g_{s}^{\prime }(r)}{\frac{\partial f}{\partial s}(r,s)},
\end{equation*}
which is positive at points $(r,s)\in \left] 0,1\right[ \times \left] 1,2\right[ $ satisfying $g_{s}(r)=1$, whence $s(r)$ is strictly increasing.

We conclude that $\dim _{\mathcal{H}}(\mathrm{supp}\left( Q_{T_{r}}\right))=\dim _{\mathcal{B}}\left( \mathrm{supp}\left( Q_{T_{r}}\right) \right) =s(r)$ is a bijection between $\left] 0,1\right[ $ and $\left] 1,2\right[ $, which completes the proof.
\end{proof}

\begin{figure}[tbh]
\begin{center}
\hspace{-2.5cm}\epsfig{file=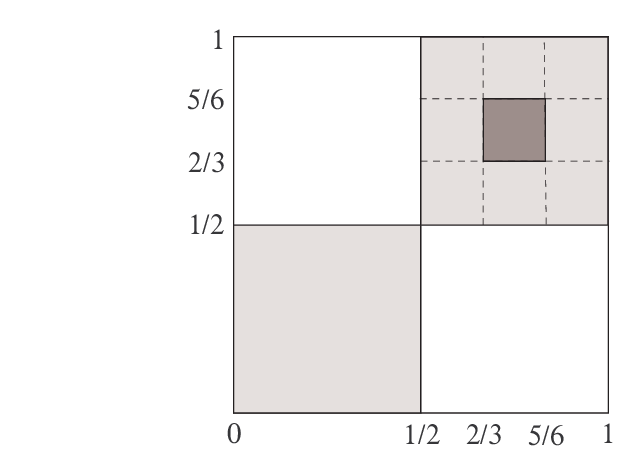,width=7cm}
\end{center}
\caption{The support of $T_{1/2}(\Pi )$ from Theorem \ref{frac1}.}
\label{fig:T4}
\end{figure}

Figure \ref{fig:T4} shows the support of the proper quasi-copula $T_{1/2}(\Pi)$ from Theorem \ref{frac1}.

In view of Theorem \ref{frac1}, a natural question arises: does there exist proper quasi-copulas which are not associated with signed measures and whose supports have dimension $2$? The answer is positive, as it is shown in \cite[Example 2.2]{Durante2016}.

\subsection{Hausdorff dimension less than 1?}

$\phantom{999}$In the literature that we have used in this work, it is assumed that the dimension of the support of a copula cannot be less than one, but we have not found any reasoning showing that this assumption is true. For this reason we want to provide a proof for that assertion. The reasoning is based on the existence of measure $\mu(C)$ induced by a copula $C$. But before providing such a proof, we need to recall the following well-known result, which can be found, for instance, in \cite{Ambrosio2008}.

\begin{lemma}[Disintegration theorem]\label{disin}
Let $\mu $\ be a measure on $\mathfrak{B}^2$, and let $\gamma$ be its projection on $\mathfrak{B}^1$. Then there exists a $\gamma$-almost everywhere uniquely determined Borel family of probability measures $\{\mu_{x}\}_{x\in [0,1]}$ on $\mathfrak{B}^1$ such that
\begin{equation*}
\mu \left( B\right) =\int_{[0,1]}\left( \int_{[0,1]}\mathbf{1}_{B_{x}}(y)\,\mathrm{d}\mu _{x}(y)\right) \mathrm{d}\gamma(x)=\int_{[0,1]}\mu_{x}(B_{x})\,\mathrm{d}\gamma(x),
\end{equation*}
for every $B\in\mathfrak{B}^2$, where $B_x=\{y\in[0,1]: (x,y)\in B\}$ and $\mathbf{1}_{B_{x}}$ denotes the indicator function of $B_{x}$.
\end{lemma}

Now we are able of proving:

\begin{theorem}\label{copless1}
For any $C\in\mathcal{C}$, we have $\dim_{\mathcal{H}}\left(\mathrm{supp}\left(\mu_C\right)\right)\ge 1$.
\end{theorem}

\begin{proof}
Lemma \ref{disin} and the fact that the projection of a set $D$ has Hausdorff dimension less than or equal to that of $D$ (see \cite[Corollary 2.4]{Falconer2014}) tell us that the marginal distributions of a (bivariate) distribution function is concentrated in sets with a Hausdorff dimension lower than $1$. Since the marginal distributions of a copula $C$ are uniform in $[0,1]$, we conclude that $\dim_{\mathcal{H}}\left(\mathrm{supp}\left(\mu_C\right)\right)\ge 1$.
\end{proof}

We can apply the same reasoning in Theorem \ref{copless1} to the case of quasi-copulas in $Q_S$; however, it cannot be applied to the quasi-copulas of the set $\mathcal{Q}_{N}$. Despite this, the result remains true for $\mathcal{Q}_{N}$, as the next result shows.

\begin{theorem}\label{thless1}
For any $Q\in \mathcal{Q}$, we have $\dim_{\mathcal{H}}\left(\mathrm{supp}\left(Q\right)\right)\ge 1$.
\end{theorem}

\begin{proof}
Suppose $\dim _{\mathcal{H}}(\mathrm{supp}(Q))<1$, then we have $\dim _{\mathcal{H}}(\pi (\mathrm{supp}(Q)))<1$, where $\pi $ is a projection of $\mathrm{supp}(Q)$ \cite[Corollary 2.4]{Falconer2014}. Therefore, there exists $x\in \left] 0,1\right[ $ such that $x\notin \pi (\mathrm{supp}(Q))$; i.e., $\left( x,y\right) \notin\mathrm{supp}(Q)$ for all $y\in \left[ 0,1\right] $. Thus, for $\left(x,y\right) $ it is possible to find an open $Q$-null rectangle $R_{\left(
x,y\right) }\in \mathcal{R}$. We cover $\left\{ x\right\} \times \left[ 0,1\right] $ with the sets $R_{\left( x,y\right) }$. With $\left\{ x\right\}\times \left[ 0,1\right] $ being a compact set, and the sets $R_{\left(x,y\right) }$ open in the topology induced in $\left[ 0,1\right] ^{2}$, we can have a finite covering of rectangles $R_{\left( x,y\right) }$. The existence of this finite covering ensures that we can find a rectangle $\left[ a,b\right] \times \left[ 0,1\right] \subset \cup R_{\left( x,y\right)}$, with $x\in \left[ a,b\right] $ and $a<b$. The inclusion $\left[ a,b
\right] \times \left[ 0,1\right] \subset \cup R_{\left( x,y\right) }$ ensures that $\left] a,b\right[ \times \left[ 0,1\right] $ is $Q$-null. On the other hand, we have $V_{Q}\left( \left[ a+\varepsilon ,b-\varepsilon\right] \times \left[ 0,1\right] \right) =b-a-2\varepsilon $ for every $\varepsilon $ in $]0,(b-a)/2[$. This contradiction leads us to the fact that the support of the quasi-copula cannot have a Hausdorff dimension less than $1$, and this completes the proof.
\end{proof}

\begin{remark}
Note that in the proof of Theorem \ref{thless1}, the fact that $Q\in\mathcal{Q}_{N}$ is not required.
\end{remark}

\section{The multidimensional case}\label{multi}

$\phantom{999}$We can obtain results similar to the previous ones for the case of multidimensional quasi-copulas. The definition of support of an $n$-quasi-copula is identical to the one that appears in the two-dimensional case given in Definition \ref{openrec} replacing rectangles with boxes of the form $I_{1}\times \cdots \times I_{n}$. By using similar ideas, we obtain that, for a fixed dimension $n$ and $s\in ]1,n]$, there exists a proper $n$-quasi-copula whose support is a self-similar set with dimension $s$.

Let $n\ge 2$ be a natural number. An $n$-{\it dimensional quasi-copula} (briefly, $n$-{\it quasi-copula}) is a function $Q\colon[0,1]^{n}\longrightarrow [0,1]$ which satisfies the following conditions:
\begin{itemize}
\item[(Q1')] For every ${\bf u}=(u_{1},u_{2},\ldots,u_{n})$ in $[0,1]^{n}$, $Q({\bf u})=0$ if at least one coordinate of {\bf u} is $0$, and $C({\bf u})=u_{k}$ whenever all coordinates of {\bf u} are $1$ except, maybe, $u_{k}$;
\item[(Q2')] $Q$ is non-decreasing in each variable;
\item[(Q3')] the {\it Lipschitz condition}: $|Q({\bf u})-Q({\bf v})|\le \sum_{i=1}^{n}|u_{i}-v_{i}|$ for all ${\bf u}, {\bf v}\in [0,1]^{n}$
\end{itemize}
(see \cite{Cuculescu2001}). The concept of $Q$-volume remains valid for $n$-{\it boxes} on $[0,1]^n$, and $Q$ is a proper $n$-quasi-copula when $Q$ is not an $n$-copula (see \cite{Durante2016book} for details).

\begin{definition}\label{matriz}
Given a natural number $n\ge 2$, $n$ divisions of $[0,1]$ of the form
\begin{equation*}
0=a_{i,0}<a_{i,1}<\cdots <a_{i,m_{i}}=1
\end{equation*}
with $m_{i}\geq 2,$ for $i=1,\ldots ,n$, and an $n$-quasi-copula $B$, we have a multidimensional matrix $T$ of dimension ${\bf m}=m_{1}\times \cdots \times m_{n} $, with $t_{i_{1},\ldots,i_{n}}=V_{B}\left( \left[ \mathbf{a}_{i_{1}-1,\ldots ,i_{n}-1},\mathbf{a}_{i_{1},\ldots ,i_{n}}\right] \right) $, where $\mathbf{a}_{i_{1}-1,\ldots ,i_{n}-1}=\left( a_{1,i_{1}-1},\ldots,a_{n,i_{n}-1}\right) $ and $\mathbf{a}_{i_{1},\ldots ,i_{n}}=\left(a_{1,i_{1}},\ldots ,a_{n,i_{n}}\right) $, and $\left[ \mathbf{a}_{i_{1}-1,\ldots ,i_{n}-1},\mathbf{a}_{i_{1},\ldots ,i_{n}}\right] $ is the set of vectors $\mathbf{u}$ such that $a_{j,i_{j}-1}\leq u_{j}\leq
a_{j,i_{j}}$. Let $\mathcal{T}_{\mathcal{Q}}$ denote the set of matrices of this type.
\end{definition}

\begin{remark}
We note that, due to \cite[Theorem 2]{Omladic2022} (see also \cite{RodUb2009}), a matrix $T=\left(t_{\mathbf{i}}\right)$ of dimension {\bf m} belongs to $\mathcal{T}_{\mathcal{Q}}$ if, and only if, it satisfies:
\begin{itemize}
\item[(a)] The sum of all entries of $T$ equals 1;
\item[(b)] $\sum_{{\bf i}\in\{j,k\}^*}t_{\bf i}>0$, where $\{j,k\}^*=\left\{{\bf i}:i_j=k\right\}$; and
\item[(c)] For any index ${\bf r}$ such that $r_j=k$,
$$\sum_{{\bf i}\in\{j,k\}^*}t_{\bf i}\ge \sum_{{\bf i}\le {\bf r}}t_{\bf i}-\sum_{{\bf i}\le {\bf r}^{(j)}}t_{\bf i}\ge 0,$$
where ${\bf r}\le {\bf i}$ denotes $r_h\le i_h$ for all $h$ and ${\bf r}^{(j)}$ is the index that coincides with ${\bf r}$ in all entries except $r_j^{(j)}=k-1$.
\end{itemize}
\end{remark}

For $\mathbf{i}=\left( i_{1},\ldots ,i_{n}\right)$, let $r_{\mathbf{i}}(\mathbf{u})$ denote the function given by
\begin{equation*}
r_{\mathbf{i}}(\mathbf{u})=\left( \max \left( \min \left( \frac{u_{1}-a_{1,i_{1}-1}}{a_{1,i_{1}}-a_{1,i_{1}-1}},1\right) ,0\right) ,\ldots
,\max \left( \min \left( \frac{u_{n}-a_{n,i_{n}-1}}{a_{n,i_{n}}-a_{n,i_{n}-1}},1\right) ,0\right) \right).
\end{equation*}

\begin{definition}\label{TT}
Given a multidimensional matrix $T$ of dimension ${\bf m}$ and an $n$-quasi-copula $Q$, we define
\begin{equation*}
T(Q)(\mathbf{u})=\sum_{\mathbf{i}}t_{\mathbf{i}}Q\left( r_{\mathbf{i}}(\mathbf{u})\right).
\end{equation*}
\end{definition}

Note that $T\in \mathcal{T}_{\mathcal{Q}}$ is not required in Definition \ref{TT}.

The following results are an adaptation of the results obtained in \cite{FerRodUb2011} to the above context.

\begin{theorem}\label{TQn}
Let $T\in \mathcal{T}_{\mathcal{Q}}$. Then, for every $n$-quasi-copula $Q$, the $T$-transformation of $Q$, $T(Q)$, is an $n$-quasi-copula. Moreover, for every $\mathbf{u}\in \left[ \mathbf{a}_{i_{1}-1,\ldots,i_{n}-1},\mathbf{a}_{i_{1},\ldots ,i_{n}}\right] $, the following equality holds:
\begin{equation*}
V_{T(Q)}([a_{1,i_{1}-1},u_{1}]\times \cdots \times \lbrack a_{n,i_{n}-1},u_{n}])=t_{\mathbf{i}}\cdot Q\left( \frac{u_{1}-a_{1,i_{1}-1}}{a_{1,i_{1}}-a_{1,i_{1}-1}},\cdots ,\frac{u_{n}-a_{n,i_{n}-1}}{a_{n,i_{n}}-a_{n,i_{n}-1}}\right) .
\end{equation*}
In particular, we have $V_{T(Q)}(\left[ \mathbf{a}_{i_{1}-1,\ldots ,i_{n}-1},\mathbf{a}_{i_{1},\ldots ,i_{n}}\right] )=t_{\mathbf{i}}$ for every $\mathbf{i}$.
\end{theorem}

From the first statement in Theorem \ref{TQn}, given $T\in \mathcal{T}_{\mathcal{Q}}$, we can define a mapping from the metric space of $n$-quasi-copulas to itself by $Q\longrightarrow T(Q)$. Without confusion, this mapping will be also denoted by $T$.

The following result ---whose proof is similar to the one given in \cite[Theorem 3.3]{FerRodUb2011} for the bivariate case, and so we omit it--- shows that, with an additional condition, the mapping $T$ is contractive, in which ${\bf i}'=\left(i_1-1,\ldots,i_n-1\right)$.

\begin{theorem}
\label{distance} Let $T=(t_{\mathbf{i}})\in \mathcal{T}_{\mathcal{Q}}$ be a matrix of order $\mathbf{m}$, and let
$$\alpha =\max \left(\sum_{{\bf r}\le{\bf i}}\left|t_{\bf r}\right|-\sum_{{\bf r}\le{\bf i}'}\left|t_{\bf r}\right|\right).$$
If $d$ denotes the sup metric, then
\begin{equation*}
d(T(Q_{1}),T(Q_{2}))\le \alpha \cdot d(Q_{1},Q_{2})  \label{contrac}
\end{equation*}
for any pair of $n$-quasi-copulas $Q_{1}$ and $Q_{2}$. As a consequence, when $\alpha <1$, we have the following:
\begin{enumerate}
\item $T$ is a contraction mapping and there is a unique $n$-quasi-copula $Q_{T}$ for which $T(Q_{T})=Q_{T}$.
\item If $T$ has a negative entry, then, for $l\in \mathbb{N}$, $T^{l}(Q)$ is a proper $n$-quasi-copula for every $n$-quasi-copula $Q$.
\item $Q_{T}$ is a proper $n$-quasi-copula if, and only if, $T$ has a negative entry.
\item If $T$ has a negative entry then $Q_{T}$ does not induce a signed measure on $[0,1]^{n}$.
\end{enumerate}
\end{theorem}

We want to note that Theorem \ref{lsupp} remains true in higher dimensions, its proof is similar with the necessary adaptation of Lemma \ref{lema} in which the marginals of the $n$-quasi-copula $Q$ will appear.

In the following, we are going to use a multidimensional matrix $T$ of dimension $4\times\overset{n}{\cdots }\times 4$. We build it in several steps. If a condition appears on the indices, it is implied in the next condition that those that appear in the previous steps are excluded.

\begin{enumerate}
\item[]Step 1. $t_{1,\ldots ,1}=1-r$,

\item[]Step 2. $t_{i_{1},\ldots ,i_{n}}=0$ if any $i_{s}=1$,

\item[]Step 3. $t_{3,\ldots ,3}=-\frac{r}{3}\left( 3^{n-1}-1+\frac{2n-1}{2n-3}\right) ^{-1}$,

\item[]Step 4. $t_{i_{1},\ldots ,i_{n}}=\frac{2n-1}{2n-3}\cdot\frac{r}{3}\left(3^{n-1}-1+\frac{2n-1}{2n-3}\right) ^{-1}$ if $i_{s}=3$ except in an index,

\item[]Step 5. $t_{i_{1},\ldots ,i_{n}}=\frac{r}{3}\left( 3^{n-1}-1+\frac{2n-1}{2n-3}\right) ^{-1}$ otherwise.
\end{enumerate}

Thus, the matrix $T$ has:
\begin{itemize}
\item[]$\rhd$ an entry equal to $1-r$,
\item[]$\rhd$ $4^{n}-3^{n}-1$ entries equal to $0$,
\item[]$\rhd$ an entry equal to $-\frac{r}{3}\left( 3^{n-1}-1+\frac{2n-1}{2n-3}\right) ^{-1}$,
\item[]$\rhd$ $2n$ entries equal to $\frac{2n-1}{2n-3}\cdot\frac{r}{3}\left(3^{n-1}-1+\frac{2n-1}{2n-3}\right) ^{-1}$ and
\item[]$\rhd$ $3^{n}-1-2n$ entries equal to $\frac{r}{3}\left( 3^{n-1}-1+\frac{2n-1}{2n-3}\right) ^{-1}$.
\end{itemize}

We observe that this matrix $T$ satisfies the condition $\alpha < 1$.

As a consequence of the results in \cite{RodUb2009} applied to $T\left( \Pi\right)$, we have that $T\left( \Pi \right) $ is an $n$-quasi-copula, whence $T\in \mathcal{T}_{\mathcal{Q}}$.

Note that we have a division of $[0,1]^{n}$ into $n$-boxes of the form $I_{1}\times \cdots \times I_{n}$ with
$$I_{s} \in \left\{ \lbrack 0,1-r],[1-r,1-2r/3],[1-2r/3,1-r/3],[1-r/3,1]\right\}.$$
Each $n$-box with equal edges (a ``hypercube") is associated with $t_{\mathbf{i}}\neq 0$; and when the box has at least two different edges, then it is associated with $t_{\mathbf{i}}=0$. That is, the support will be self-similar. Studying the equation $\left( 1-r\right) ^{s}+3^{n}\left( \frac{r}{3}\right) ^{s}=1$ in a similar way as in the proof of Theorem \ref{frac1} for $n=2$, we obtain the following result:

\begin{theorem}
For all $s\in ]1,n[$ there exists a proper $n$-quasi-copula whose support is self-similar and its dimension is $s$.
\end{theorem}

Finally, we observe that the case $s=n$ is obtained with the matrix of dimension $3\times \overset{n}{\cdots }\times 3$ given by:
\begin{eqnarray*}
t_{2,\ldots ,2}\!\!\!&=&\!\!\!-\frac{1}{3}\left( 3^{n-1}-1+\frac{2n-1}{2n-3}\right)^{-1},\\
t_{i_{1},\ldots ,i_{n}}\!\!\!&=&\!\!\!\left\{\begin{array}{ll}
\frac{2n-1}{2n-3}\cdot\frac{1}{3}\left( 3^{n-1}-1+\frac{2n-1}{2n-3}\right) ^{-1}, & \text{if $i_{s}=2$ except in an index},\\ \noalign{\medskip}
\frac{1}{3}\left( 3^{n-1}-1+\frac{2n-1}{2n-3}\right) ^{-1},&  \text{otherwise}.
\end{array} \right.
\end{eqnarray*}

\section{Conclusions}\label{conclus}

$\phantom{999}$In this paper, we have introduced the concept of support of a quasi-copula, both of $\mathcal{Q}_{N}$ and $\mathcal{Q}_{S}$, in order to prove that for every $s\in ]1,2[$ there is a (proper) quasi-copula in ${\mathcal{Q}}_{N}$ whose support is a self-similar fractal set with Hausdorff dimension ---and box-counting dimension--- equal to $s$. It remains an open problem to establish whether the result is still true when we substitute ${\mathcal{Q}}_{N}$ by ${\mathcal{Q}}_{S}$. We have also addressed the multidimensional case, generalizing results from the bidimensional case.
\bigskip
\smallskip

\noindent{\large \textbf{Acknowledgements}}\medskip

$\phantom{999}$ The second and third authors also acknowledge the support of the program FEDER-Andaluc\'{\i}a 2014-2020 (Spain) under research project UAL2020-AGR-B1783 and the support of the project PID2021-122657OB-I00 by the Ministerio de Ciencia e Innovaci\'{o}n (Spain).


\smallskip

\begin{thebibliography}{99}
\bibitem{Alsina93} C. Alsina, R.B. Nelsen, B. Schweizer, On the characterization of a class of binary operations on distribution functions, Statist. Probab. Lett. 17 (1993) 85--89.

\bibitem{Ambrosio2008} L. Ambrosio, N. Gigli, G. Savar\'e, Gradient Flows in Metric Spaces and in the Space of Probability Measures, 2nd ed., Lectures in Mathematics ETH Z\"{u}rich, Birkh\"{a}user Verlag, Basel, 2008.

\bibitem{Arias2020} J.J. Arias-Garc\'ia, R. Mesiar, B. De Baets, A hitchhiker's guide to quasi-copulas, Fuzzy Sets Syst. 393 (2020) 1--28.

\bibitem{Barnsley1989} M.F. Barnsley, Lecture Notes on Iterated Function Systems, in: R.L. Devaney, L. Keen, L. (Eds.), Chaos and Fractals: The Mathematics Behind the Computer Graphics, Proceedings of Symposia on Applied Mathematics, vol. 39, American Mathematical Society, Providence, 1989, pp. 127--144.

\bibitem{Beliakov2007} G. Beliakov, A. Pradera, T. Calvo, Aggregation Functions: A Guide for Practitioners, Studies in Fuzziness and Soft Computing, vol.221, Springer, Berlin, 2007.

\bibitem{Cuculescu2001} I. Cuculescu, R. Theodorescu, Copulas: diagonals and tracks, Rev. Roum. Math. Pures Appl. 46 (2001) 731--742.

\bibitem{deAmo2012} E. de Amo, M. D\'{\i}az Carrillo, J. Fern\'andez S\'anchez, Characterization of all copulas associated with non-continuous random variables, Fuzzy Sets Syst. 191 (2012) 103--112.

\bibitem{deAmo2012b} E. de Amo, M. D\'iaz Carrillo, J. Fern\'andez-S\'anchez, Copulas and associated fractal sets, J. Math. Anal. Appl. 386 (2012), 528--541.

\bibitem{Durante2016book} F. Durante, C. Sempi, Principles of Copula Theory, Chapman $\&$ Hall/CRC, Boca Raton, 2016.

\bibitem{Durante2016} F. Durante, J. Fern\'andez-S\'anchez, W. Trutschnig, Baire category results for quasi-copulas, Depend. Model. 4 (2016), 215--223.

\bibitem{Edgard2008} G. Edgar, Measure, Topology, and Fractal Geometry, 2nd ed., Springer-Verlag, New York, 2008.

\bibitem{Falconer2014} K. Falconer, Fractal Geometry: Mathematical Foundations and Applications, 3rd ed., Wiley, Hoboken, 2014.

\bibitem{FerUb2014} J. Fern\'andez-S\'anchez, M. \'Ubeda-Flores, A note on quasi-copulas and signed measures, Fuzzy Sets Syst. 234 (2014) 109--112.

\bibitem{FerNelUb2011} J. Fern\'andez-S\'anchez, R.B. Nelsen, M. \'Ubeda-Flores, Multivariate copulas, quasi-copulas, and lattices, Statist. Probab. Lett. 81 (2011) 1365--1369.

\bibitem{FerRodUb2011} J. Fern\'{a}ndez-S\'{a}nchez, J.A. Rodr\'{\i}guez-Lallena, M. \'{U}beda-Flores, Bivariate quasi-copulas and doubly stochastic signed measures, Fuzzy Sets Syst. 168 (2011) 81--88.

\bibitem{Fredricks2005} G.A. Fredricks, R.B. Nelsen, J.A. Rodr\'iguez-Lallena, Copulas with fractal supports, Insurance Math. Econom. 37 (2005) 42--48.

\bibitem{Genest99} C. Genest, J.J. Quesada-Molina, J.A. Rodr\'{\i}guez-Lallena, C. Sempi, A characterization of quasi-copulas, J. Multivariate Anal. 69 (1999) 193--205.

\bibitem{Grabish2009} M. Grabisch, J.L. Marichal, R. Mesiar, E. Pap, Aggregation Functions, Encyclopedia of Mathematics and its Applications, vol. 127. Cambridge University Press, Cambridge, 2009.

\bibitem{Halmos1974} P.R. Halmos, Measure Theory, Springer, New York, 1974.

\bibitem{Hutchinson1981} J.E. Hutchinson, Self-Similar Sets, Indiana Univ. Math. J. 30 (1981) 713--747.

\bibitem{Mandelbrot1977} B.B. Mandelbrot, Fractals: Form, Chance and Dimension, W.H. Freeman and Co., New York, 1977.

\bibitem{Nelsen2006} R.B. Nelsen, An Introduction to Copulas, 2nd ed., Springer, New York, 2006.

\bibitem{Nelsen2010} R.B. Nelsen, J.J. Quesada-Molina, J.A. Rodr\'{\i}guez-Lallena, M. \'{U}beda-Flores, Quasi-copulas and signed measures, Fuzzy Sets Syst. 161 (2010) 2328--2336.
    
\bibitem{Omladic2022}M. Omladi$\check{\rm c}$, N. Stopar, Multivariate imprecise Sklar type theorems, Fuzzy Sets Syst. 428 (2022) 80--101.

\bibitem{RodUb2009} J.A. Rodr\'{\i}guez-Lallena, M. \'{U}beda-Flores, M, Some new characterizations and properties of quasi-copulas, Fuzzy Sets Syst. 160 (2009) 717--725.

\bibitem{Sempi2017} C. Sempi, Quasi-copulas: a brief survey, in: M. \'Ubeda-Flores, E. de Amo Artero, F. Durante, J. Fern\'andez-S\'anchez (Eds.), Copulas and Dependence Models with Applications, Springer, Cham, 2017, pp. 203--224.

\bibitem{Sklar59} A. Sklar, Fonctions de r\'epartition $\grave{\mathrm{a}}$ $n$ dimensions et leurs marges, Publ. Inst. Statist. Univ. Paris 8 (1959) 229--231.

\bibitem{Trutschnig2012} W. Trutschnig, J. Fern\'{a}ndez-S\'{a}nchez, Idempotent and multivariate copulas with fractal support. J. Statist. Plann. Inference 142 (2012) 3086--3096.
\end{thebibliography}
\end{document}